\documentclass{scrartcl}
\usepackage[utf8]{inputenc}
\usepackage[T1]{fontenc}
\usepackage{lmodern}

\usepackage[defthms,stable]{arw-131203}

\newcommand{\Ast}{\mathcal{A}}

\newcommand*{\et}[1]{\mathbf{e}_{#1}}

\newcommand{\rhoh}[0]{\hat{\rho}}

\newvarcmd{\Indic}{\Ind\sizel[#2\sizer]}

\newcommand{\Ysu}{\overline{Y}}

\newcommand{\sigzo}{\sigma^{[0,1]}}
\newcommand{\Uo}{U_1}
\newcommand{\uo}{u_1}

\hypersetup{
  pdfauthor={A. E. Kyprianou, A. R. Watson},
  pdftitle={Potentials of stable processes}
}

\title{Potentials of stable processes}

\author{%
  A.\ E.\ Kyprianou%
  \footnote{%
    University of Bath, UK. \email{a.kyprianou@bath.ac.uk}%
  }
  \and A.\ R.\ Watson%
  \footnote{%
    CIMAT, Mexico. \email{alexander.watson@cimat.mx}.
    This work was begun while the second author was a Ph.D.\ student at the
    University of Bath, and completed during a position at CIMAT.%
  }%
}

\begin{document}

\maketitle

\begin{abstract}
\noindent \textbf{Abstract.}
For a stable process, we
give an explicit formula for the potential measure of the
process killed outside a bounded interval
and the joint law
of the overshoot, undershoot and undershoot from the maximum
at exit from a bounded interval.
We obtain the equivalent quantities for a stable process
reflected in its infimum.
The results are obtained by exploiting a
simple connection with the Lamperti representation
and exit problems of stable processes.

{\small

}
\end{abstract}

\section{Introduction and results}

For a L\'evy process $X$, the measure
\[ U^A(x, \dd y) = \stE_x \int_0^\infty
  \Indic{X_t \in \dd y} \Indic{\forall s \le t : X_s \in A}\, \dd t,
\]
called the
\define{potential} (or \define{resolvent})
\define{measure of $X$ killed outside $A$},
is a quantity of great interest, and is related to exit problems.

The main cases where the potential measure can be
computed explicitly are as follows.
If $X$ is a L\'evy process with known Wiener--Hopf factors,
it can be obtained when $A$ is half-line or $\RR$;
see \cite[Theorem VI.20]{Ber-Levy}.
When $X$ is a totally asymmetric L\'evy process with known
scale functions, it can be obtained for $A$ a bounded interval,
a half-line or $\RR$; see \cite[Section 8.4]{Kyp}.
Finally, \cite{Bau-reflected} details a technique to obtain
a potential measure for a reflected L\'evy process killed outside a
bounded interval from the same quantity for the
unreflected process.

In this note, we consider the case where $X$ is a stable process
and $A$ is a bounded interval. We compute the measure $U^{[0,1]}$,
from which $U^A$ may be obtained for any bounded interval $A$
via spatial homogeneity and scaling;
and from this we compute the joint law at first exit
of $[0,1]$ of the overshoot,
undershoot and undershoot from the maximum.
Furthermore, we give the potential measure and triple law
also for the process reflected in its infimum.

The potential measure has been already been computed when $X$ is
symmetric; see \citet[Corollary 4]{BGR-st-hit} and references therein,
as well as \citet{Bau-reflected}.
We extend these results to asymmetric stable processes
with jumps on both sides. The essential
observation is that a potential for $X$ with killing outside a bounded
interval may be converted into a potential for the
\emph{Lamperti transform of $X$}, say $\xi$,
 with killing outside a half-line.
To compute this potential in a half-line, it is enough
to know the killing rate of $\xi$ and the solution
of certain exit problems for $X$.
The results for the reflected process are obtained
via the work of \citet{Bau-reflected}.


We now give our results. Some facts we will rely on
are summarised in \autoref{s:prelim}, and proofs are 
given in \autoref{s:proofs}.

\medskip\noindent
We work with the (strictly) stable process $X$ with scaling parameter $\alpha$ and
positivity parameter $\rho$, which is defined as follows. For $(\alpha,\rho)$
in the set
\begin{eqnarr*}
  \Ast
  &=&
  \{ (\alpha,\rho) : \alpha \in (0,1), \, \rho \in (0,1) \}
  \cup \{ (\alpha,\rho) = (1,1/2) \}\\
  &&
  {} \cup \{ (\alpha,\rho) : \alpha \in (1,2), \, \rho \in (1-1/\alpha, 1/\alpha) \},
\end{eqnarr*}
let $X$, with probability laws $(\stP_x)_{x \in \RR}$,
be the L\'evy process with characteristic exponent
\[ \CE(\theta)
  = \begin{cases}
    c\abs{\theta}^\alpha (1  - \iu\beta\tan\frac{\pi\alpha}{2}\sgn\theta)
    & \alpha \in (0,2) \setminus \{1\}, \\
    c\abs{\theta} & \alpha = 1,
  \end{cases}
  \wh \theta \in \RR,
\]
where $c = \cos(\pi\alpha(\rho-1/2))$ and
$\beta = \tan(\pi\alpha(\rho-1/2))/\tan(\pi\alpha/2)$.
This L\'evy process has absolutely continuous L\'evy measure with density
\newcommand{\cp}{\frac{\Gamma(\alpha+1)}{\Gamma(\alpha\rho)\Gamma(1-\alpha\rho)}}
\newcommand{\cm}{\frac{\Gamma(\alpha+1)}{\Gamma(\alpha\rhoh)\Gamma(1-\alpha\rhoh)}}
\[ c_+ x^{-(\alpha+1)} \Indic{x > 0}
  + c_- \abs{x}^{-(\alpha+1)}\Indic{x < 0},
  \wh x \in \RR , \]
where
\[ c_+ = \cp, \qquad c_- = \cm \]
and $\rhoh = 1-\rho$.

The parameter set $\Ast$ and the characteristic exponent $\CE$ represent,
up a multiplicative constant in $\Psi$, all (strictly) stable processes
which jump in both directions, except for
Brownian motion and the symmetric Cauchy processes with non-zero drift.
The normalisation is the same as that in \cite{KP-HG}, and when $X$
is symmetric, that is when $\rho = 1/2$, the normalisation
agrees with that of
\cite{BGR-st-hit}.
We remark that the quantities we are interested in
can also be derived in cases of one-sided jumps:
either $X$ is a subordinator, in which case the results are
trivial, or $X$ is a spectrally one-sided L\'evy process, in which case
the potentials in question may be assembled using the theory of scale
functions; see \cite[Theorem 8.7 and Exercise 8.2]{Kyp}.

The choice $\alpha$ and $\rho$ as parameters is explained as follows.
$X$ satisfies the \define{$\alpha$-scaling property}, that
\begin{equation}\label{e:scaling}
  \text{under } \stP_x \text{, the law of }
  (cX_{t c^{-\alpha}})_{t \ge 0} \text{ is } \stP_{cx} ,
\end{equation}
for all $x \in \RR, \, c > 0$. The second parameter
satisfies $\rho = \stP_0(X_t > 0)$.

\medskip\noindent
Having defined the stable process, we proceed to our results. Let
\[ \sigzo = \inf \{ t \ge 0 : X_t \notin [0,1] \}, \]
and define the killed potential measure and potential density
\[ \Uo(x, \dd y)
  := U^{[0,1]}(x, \dd y)
  = \stE_x \int_0^{\sigzo} \Indic[var]{X_t \in \dd y} \, \dd t
  = \uo(x,y) \, \dd y, \]
provided the density $\uo$ exists.

\begin{thrm}\label{X (0,1)}
For $0 < x, y < 1$,
\[ \uo(x,y)
  = \begin{cases}
    \dfrac{1}{\Gamma(\alpha\rho)\Gamma(\alpha\rhoh)}
    (x-y)^{\alpha-1}
    \dint_0^{\frac{y(1-x)}{x-y}}
    s^{\alpha\rho-1}
    (s+1)^{\alpha\rhoh-1} \, \dd s,
    & y < x, \\
    \dfrac{1}{\Gamma(\alpha\rho)\Gamma(\alpha\rhoh)}
    (y-x)^{\alpha-1}
    \dint_0^{\frac{x(1-y)}{y-x}}
    s^{\alpha\rhoh-1}
    (s+1)^{\alpha\rho-1} \, \dd s,
    & x < y.
  \end{cases}
\]
\end{thrm}
When $X$ is symmetric, this reduces,
by spatial homogeneity and scaling of $X$,
and substituting in the integral,
to \cite[Corollary 4]{BGR-st-hit}.

\newcommand{\tbp}{\tau_1^+}
\newcommand{\tzm}{\tau_0^-}
With very little extra work, \autoref{X (0,1)} yields an
apparently stronger result. Let
\[ \tzm = \inf\{ t \ge 0: X_t < 0\} ;
  \qquad \qquad \Xsu_t = \sup_{s \le t} X_s, \wh t \ge 0,
\]
and write
\[
  \stE_x \int_0^{\tzm}
  \Indic{X_t \in \dd y, \, \Xsu_t \in \dd z} \, \dd t
  = u(x,y,z) \, \dd y \, \dd z,
\]
if the right-hand side exists. Then we have the following.
\begin{cor}
  For $x > 0$, $y \in [0,z)$, $z > x$,
  \begin{equation}\label{e:X double pot}
    u(x,y,z)
    =
    \frac{1}{\Gamma(\alpha\rho)\Gamma(\alpha\rhoh)}
    x^{\alpha\rhoh}
    y^{\alpha\rho}
    \frac{
      (z-x)^{\alpha\rho-1}
      (z-y)^{\alpha\rhoh-1}
    }{
      z^{\alpha}
    }
    \, \dd y \, \dd z.
  \end{equation}
  \begin{proof}
    Rescaling, we obtain
    \[
      \stE_x \int_0^{\tzm} \Indic{X_t \in \dd y, \, \Xsu_t \le z} \, \dd t
      = z^{\alpha-1} u_1(x/z,y/z),
    \]
    and the density is found by differentiating the right-hand side in $z$.
  \end{proof}
\end{cor}
From this density, one may recover the following hitting distribution,
which originally appeared in \citet[Corollary 15]{KPR-n-tuple}. Let
\[
  \tbp = \inf\{ t \ge 0: X_t > 1\}. \]
\begin{cor}\label{c:X hit}
  For $u \in [0,1-x)$, $v \in (u,1]$, $y \ge 0$,
  \begin{multline}\label{e:X hit}
    \stP_x(1-\Xsu_{\tbp-} \in \dd u, \, 1-X_{\tbp-}\in \dd v, \,
    X_{\tbp} - 1 \in \dd y, \, \tbp<\tzm) \\
    =
    \frac{\Gamma(\alpha+1)}{\Gamma(\alpha\rhoh)\Gamma(1-\alpha\rho)}
    \frac{
      x^{\alpha\rhoh}
      (1-v)^{\alpha\rho}
      (1-u-x)^{\alpha\rho-1}
      (v-u)^{\alpha\rhoh-1}
    }{
      (1-u)^{\alpha}
      (v+y)^{\alpha+1}
    }
    \, \dd u \, \dd v \, \dd y.
  \end{multline}
  \begin{proof}
    Following the proof of \cite[Proposition III.2]{Ber-Levy},
    one may show that the left-hand side of \eqref{e:X hit} is equal to
    $u(x,1-v,1-u)\pi(v+y)$,
    where $\pi$ is the L\'evy density of $X$.
  \end{proof}
\end{cor}

\begin{rem}
  The proof of \autoref{c:X hit} suggests an alternative derivation of
  \autoref{X (0,1)}. Since the identity \eqref{e:X hit} is
  already known, one may deduce $u(x,y,z)$ from it by following
  the proof backwards. The potential $u_1(x,y)$ without $\Xsu$ may then be
  obtained via integration.
  However, in \autoref{s:proofs} we offer instead a self-contained
  proof based on well-known hitting distributions for the stable
  process.
\end{rem}

\medskip\noindent
Now let $Y$ denote the stable process $X$ reflected in its infimum,
that is,
\[ Y_t = X_t - \Xin_t, \wh t \ge 0, \]
where $\Xin_t = \inf\{X_s, 0 \le s \le t\} \meet 0$ for $t \ge 0$.
$Y$ is a self-similar Markov process.

Let $T_1^+ = \inf\{ t > 0: Y_t > 1\}$ denote the first passage time
of $Y$ above the level $1$, and define
\[ R_1(x, \dd y)
  = \stE_x \int_0^{T_1^+} \Indic[var]{Y_t \in \dd y} \, \dd t
  = r_1(x,y) \, \dd y , \]
provided that the density on the right-hand side exists.
We may then use the results of \citet{Bau-reflected} to find $r_1$.
Note that, as $Y$ is self-similar, $R_1$ suffices to deduce
the potential of $Y$ killed at first passage above any level.

\begin{thrm}\label{Y (0,1)}
For $0 < y < 1$,
\[ r_1(0,y)
  = \frac{1}{\Gamma(\alpha)} y^{\alpha\rho-1} (1-y)^{\alpha\rhoh} . \]
Hence, for $0 < x,\,y < 1$,
\[
  r_1(x,y)
  =
  \left\{
  \begin{IEEEeqnarraybox}[][c]{ll?l}
    \IEEEstrut
    \dfrac{1}{\Gamma(\alpha\rho)\Gamma(\alpha\rhoh)}
    \biggl[
    &(x-y)^{\alpha-1}
    \dint_0^{\frac{y(1-x)}{x-y}}
    s^{\alpha\rho-1}
    (s+1)^{\alpha\rhoh-1} \, \dd s
    \biggr. & \\
    &
    \biggl.
    {} + y^{\alpha\rho-1}
    (1-y)^{\alpha\rhoh}
    \dint_0^{1-x}
    t^{\alpha\rho-1}
    (1-t)^{\alpha\rhoh-1} \, \dd t
    \biggr],
    & y < x, \\
    \dfrac{1}{\Gamma(\alpha\rho)\Gamma(\alpha\rhoh)}
    \biggl[
    &(y-x)^{\alpha-1}
    \dint_0^{\frac{x(1-y)}{y-x}}
    s^{\alpha\rhoh-1}
    (s+1)^{\alpha\rho-1} \, \dd s
    \biggr. &\\
    & \biggl.
    {} + y^{\alpha\rho-1}
    (1-y)^{\alpha\rhoh}
    \dint_0^{1-x}
    t^{\alpha\rho-1}
    (1-t)^{\alpha\rhoh-1} \, \dd t
    \biggr],
    & x < y.
    \IEEEstrut
  \end{IEEEeqnarraybox}
  \right.
\]
\end{thrm}
Writing
\[
  \stE_x\int_0^\infty \Indic{Y_t \in \dd y, \, \Ysu_t \in \dd z} \, \dd t
  = r(x,y,z) \, \dd y \, \dd z,
\]
where $\Ysu_t$ is the supremum of $Y$ up to time $t$,
we obtain the following corollary, much as we had for $X$.
\begin{cor}
  For $y \in (0,z)$, $z \ge 0$,
  \[
    r(0,y,z)
    =
    \frac{\alpha\rhoh}{\Gamma(\alpha)}
    y^{\alpha\rho-1}
    (z-y)^{\alpha\rhoh-1},
  \]
  and for $x > 0$, $y \in (0,z)$, $z \ge x$,
  \begin{IEEEeqnarray*}{ll}
    r(x,y,z)
    =
    \dfrac{1}{\Gamma(\alpha\rho)\Gamma(\alpha\rhoh)}
    y^{\alpha\rho-1}
    (z-y)^{\alpha\rhoh-1}
    \biggl[
    &
    x^{\alpha\rhoh}
    (z-x)^{\alpha\rho-1}
    z^{1-\alpha}
    \\
    &{} +
    \alpha\rhoh
    \int_0^{1-\frac{x}{z}}
    t^{\alpha\rho-1}
    (1-t)^{\alpha\rhoh-1}
    \, \dd t
    \biggr]
  \end{IEEEeqnarray*}
\end{cor}

We also have the following corollary, which is the analogue of
\autoref{c:X hit}.
\begin{cor}\label{c:Y hit}
  For $u \in (0,1]$, $v \in (u,1)$, $y \ge 0$,
  \begin{multline*}
    \stP_0( 1-\Ysu_{T_1^+} \in \dd u,
    \, 1-Y_{T_1^+} \in \dd v,
    \, Y_{T_1^+} -1 \in \dd y)\\
    = \frac{\alpha \cdot \alpha\rhoh}{\Gamma(\alpha\rho)\Gamma(1-\alpha\rho)}
    \frac{(1-v)^{\alpha\rho-1}(v-u)^{\alpha\rhoh-1}}{(v+y)^{\alpha+1}}
    \dd u\, \dd v\, \dd y 
    ,
  \end{multline*}
  and for $x \ge 0$, $u \in [0,1-x)$, $v \in (u,1)$, $y \ge 0$,
  \begin{IEEEeqnarray*}{ll}
    \IEEEeqnarraymulticol{2}{l}{
      \stP_x( 1-\Ysu_{T_1^+} \in \dd u,
      \, 1-Y_{T_1^+} \in \dd v,
      \, Y_{T_1^+} -1 \in \dd y)}\\
    \quad {} = {} & \frac{\Gamma(\alpha+1)}{\Gamma(\alpha\rhoh)\Gamma(1-\alpha\rho)}
    \frac{
      (1-v)^{\alpha\rho-1}
      (v-u)^{\alpha\rhoh-1}
    }{
      (v+y)^{\alpha+1}
    } \\
    & {} \times
    \biggl[    
    x^{\alpha\rhoh}
    (1-u-x)^{\alpha\rho-1}
    (1-u)^{1-\alpha}
     + 
    \alpha\rhoh
    \int_0^{1-\frac{x}{1-u}}
    t^{\alpha\rho-1}
    (1-t)^{\alpha\rhoh-1}
    \, \dd t
    \biggr]
    \dd u\, \dd v\, \dd y 
    .
  \end{IEEEeqnarray*}
\end{cor}
The marginal in $\dd v\,\dd y$ appears in
\citet[Corollary 3.5]{Bau-reflected} for the case where $X$ is
symmetric and $x=0$. The marginal in
$\dd y$ is given in \citet{Kyp-reflected} for the process
reflected in the supremum; this corresponds to
swapping $\rho$ and $\rhoh$. However, unless $x=0$,
it appears to be difficult to
integrate in \autoref{c:Y hit} and obtain the expression found
in \cite{Kyp-reflected}.

Finally, one may integrate in \autoref{Y (0,1)} and obtain the
expected first passage time for the reflected process.
\begin{cor}
For $x \ge 0$,
\[ \stE_x[T_1^+]
  = \frac{1}{\Gamma(\alpha+1)}
  \biggl[
  x^{\alpha\rhoh}
  (1-x)^{\alpha\rho}
  + \alpha\rhoh
  \int_0^{1-x}
  t^{\alpha\rho-1}
  (1-t)^{\alpha\rhoh-1}
  \, \dd t
  \biggr].
\]
In particular,
\[ \stE_0[T_1^+]
  = \frac{1}{\Gamma(\alpha)}
  \frac{\Gamma(\alpha\rho)\Gamma(\alpha\rhoh+1)}{\Gamma(\alpha+1)} .
\]
\end{cor}

\section{The Lamperti representation}
\label{s:prelim}

\medskip\noindent
We will calculate potentials related to $X$ by appealing to the Lamperti
transform \cite{Lam-ssLT,VA-Ito}. Recall that a process $Y$ with probability
measures $(\stP_x)_{x > 0}$ is a
\define{positive self-similar Markov process (pssMp)}
if it is a standard Markov process (in the sense of \cite{BG-mppt})
with state space $[0,\infty)$
which has zero as an absorbing state and
satisfies the scaling property:
\begin{equation*}
  \text{under } \stP_x \text{, the law of }
  (cY_{t c^{-\alpha}})_{t \ge 0} \text{ is } \stP_{cx} ,
\end{equation*}
for all $x,\,c > 0$.
  
The Lamperti transform gives a correspondence between pssMps and killed
L\'evy processes, as follows. Let
$S(t) = \int_0^t (Y_u)^{-\alpha} \, \dd u$;
this process is continuous and strictly increasing until $Y$ reaches zero.
Let $T$ be its inverse. Then, the process
\[ \xi_s = \log Y_{T(s)}, \wh s \ge 0 \]
is a L\'evy process, possibly killed at an independent exponential time,
and termed the \define{Lamperti transform} of $Y$. Note that
$\xi_0 = \log x$ when $Y_0 = x$, and one may easily see from the
definition of $S$ that $e^{\alpha \xi_{S(t)}} \, \dd S(t) = \dd t$.

\medskip\noindent
A simple example of the Lamperti transform in action is given by
considering the process $X$. Define
\[ \tau_0^- = \inf\{ t \ge 0 : X_t < 0 \}, \]
and let
\[ \stPk_x( X_t \in \cdot)
  = \stP_x(X_t \in \cdot, \, t < \tau_0^-), \wh t \ge 0, \; x > 0. \]
The process $X$ with laws $(\stPk_x)_{x > 0}$
is a pssMp. \citet{CC-LS} gives explicitly the generator of its
Lamperti transform, whose laws we denote
$(\LevPk_y)_{y \in \RR}$, finding in particular that it is killed
at rate
\begin{equation}\label{e:q}
  q
  := c_-/\alpha
  = \frac{\Gamma(\alpha)}{\Gamma(\alpha\rhoh)\Gamma(1-\alpha\rhoh)}
  .
\end{equation}

\section{Proofs}
\label{s:proofs}

To avoid the proliferation of
symbols, we generally distinguish processes only by the measures
associated with them; the exception is that self-similar processes
will be distinguished from
processes obtained by Lamperti transform.
Thus, the time
\[ \tau_1^+ = \inf\{ t \ge 0 : X_t > 0 \} \]
always refers to the canonical process of the measure it appears
under, and will be used for self-similar processes; and
\[ S_0^+ = \inf\{ s \ge 0: \xi_s > 0 \},
  \quad \text{and} \quad
  S_0^- = \inf\{ s \ge 0 : \xi_s < 0 \} \]
will likewise be used for processes obtained by Lamperti transform.

\begin{proofof}{\autoref{X (0,1)}}
Our proof makes use of the pssMp $(X,\stPk)$
and its Lamperti transform $(\xi,\LevPk)$,
both defined in \autoref{s:prelim}.
Let $0 < x,\,y < 1$.
Then
\begin{eqnarr*}
  \Uo(x, \dd y)
  &=& \stE_{x} \int_0^{\sigzo} \Indic[var]{X_t \in \dd y} \, \dd t \\
  &=& \stEk_{x} \int_0^{\tau_1^+} \Indic[var]{X_t \in \dd y} \, \dd t,
\end{eqnarr*}
using nothing more than the definition of $(X,\stPk)$.
We now use the Lamperti representation to relate this to $(\xi,\LevPk)$.
This process is killed at the
rate $q$ given in \eqref{e:q}, and so it may be represented as
an unkilled L\'evy process $(\xi,\LevP)$ which is sent to some cemetery
state at the independent exponental time $\et{q}$.
We now make the following calculation, in which
$S$ and $T$ are the time-changes used in the Lamperti transform as described
in \autoref{s:prelim},
\begin{eqnarr*}
  \Uo(x, \dd y)
  &=& \LevEk_{\log(x)} \int_0^{T(S_0^+)}
    \Indic[var]{e^{\xi_{S(t)}} \in \dd y}
    e^{\alpha \xi_{S(t)}} \,
    \dd S(t) \\
  &=& y^{\alpha}
    \LevE_{\log(x)} \int_0^{S_0^+}
    \Indic[var]{e^{\xi_s} \in \dd y}
    \Indic[var]{\et{q} > s} \, \dd s \\
  &=& y^{\alpha}
    \LevEh_{\log(1/x)}
    \int_0^{S_0^-}
    \Indic[var]{\xi_s \in \log(1/\dd y)}
    e^{-qs}
    \, \dd s ,
\end{eqnarr*}
where $\LevEh$ refers to the dual L\'evy process.
Examining the proof of Theorem VI.20 in \citet{Ber-Levy}
reveals that, for any $a > 0$,
\[ \LevEh_a \int_0^{S_0^-}
  \Indic[var]{\xi_s \in \cdot} \, e^{-qs} \, \dd s
  = \frac{1}{q} \int_{[0,\infty)} \LevPh_0( \xis_{\et{q}} \in \dd w)
  \int_{[0,a]} \LevPh_0(-\xii_{\et{q}} \in \dd z)
  \Indic[var]{a + w - z \in \cdot} ,
\]
where for each $t \ge 0$, $\xis_{t} = \sup\{\xi_s : s \le t\}$
and $\xii_{t} = \inf\{\xi_s : s \le t\}$.
Then, provided that the measures $\LevPh_0(\xis_{\et{q}} \in \cdot)$
and $\LevPh_0(\xii_{\et{q}} \in \cdot)$ possess respective densities
$g_S$ and $g_I$ (as we will shortly see they do), it follows that
for $a > 0$,
\[ \LevEh_a \int_0^{S_0^-}
  \Indic[var]{\xi_s \in \dd v} \, e^{-qs} \, \dd s
  = \frac{\dd v}{q} \int_{(a-v) \join 0}^{a} \dd z \, g_I(-z) g_S(v-a+z) . \]

We may apply this result to our potential measure $\Uo$ in
order to find its density, giving
\begin{equation}\label{e:inter1}
  \uo(x,y)
  = \frac{1}{q}
  y^{\alpha-1}
  \int_{\frac{y}{x} \join 1}^{\frac{1}{x}}
  t^{-1}
  g_I(\log t^{-1})
  g_S( \log(tx/y)) \, \dd t .
\end{equation}

It remains to determine the densities $g_S$ and $g_I$ of
the measures $\LevPh_0(\xis_{\et{q}} \in \cdot)$
and $\LevPh_0(\xii_{\et{q}} \in \cdot)$.
These can be related
to functionals of $X$ by the Lamperti transform:
\begin{equation}\label{e:laws}
  \begin{IEEEeqnarraybox}[][c]{rCcCl}
    \LevPh_0(\xis_{\et{q}} \in \cdot)
    &=& \LevP_0(- \xii_{\et{q}} \in \cdot)
    &=& \stP_1(-\log \Xin_{{\tau_0^-}-} \in \cdot) \\
    \LevPh_0(\xii_{\et{q}}\in \cdot)
    &=& \LevP_0(- \xis_{\et{q}} \in \cdot)
    &=& \stP_1(-\log \Xsu_{\tau_0^-} \in \cdot),
  \end{IEEEeqnarraybox}
\end{equation}
where $\Xin$ and $\Xsu$ are defined in the obvious manner.

The laws of the rightmost random variables in \eqref{e:laws}
are available explicitly,
as we now show. For the law of $\Xin_{\tau_0^-{-}}$, we transform it
into an overshoot problem and make
use of Example 7 in \citet{DK-os}, as follows.
We omit the calculation of the integral, which uses \cite[8.380.1]{GR}.
\begin{eqnarr*}
  \stP_1(\Xin_{{\tau_0^-}-} \in \dd y)
  &=& \stPh_0(1-\Xsu_{{\tau_1^+}-} \in \dd y) \\
  &=& K \int_y^\infty \dd v \,
    \int_0^\infty \dd u \,
    (v-y)^{\alpha\rho-1} (v+u)^{-(\alpha+1)}
    (1-y)^{\alpha\rhoh-1} \, \dd y \\
  &=& \frac{\sin(\pi\alpha\rhoh)}{\pi}
    y^{-\alpha\rhoh}
    (1-y)^{\alpha\rhoh-1} \, \dd y .
  \yesnumber\label{e:law1}
\end{eqnarr*}

For the law of $\Xsu_{\tau_0^-}$, consider the following calculation.
\[
  \stP_1( \Xsu_{\tau_0^-} \ge y)
  = \stP_1( \tau_y^+ < \tau_0^-) = \stP_{1/y}(\tau_1^+ < \tau_0^-) .
\]
This final quantity depends on the solution of the two-sided
exit problem for the stable process;
it is computed in \citet{Rog-hit}, where it is denoted
$f_1(1/y,\infty)$.
Note that \cite{Rog-hit} contains a typographical error:
in Lemma 3 of that work and the discussion after it,
the roles of $q$ (which is $\rho$ in our notation)
and $1-q$ should be swapped.
In the corrected form, we have
\begin{eqnarr*}
  \stP_1( \Xsu_{\tau_0^-} \ge y)
  &=& \frac{\Gamma(\alpha)}{\Gamma(\alpha\rho)\Gamma(\alpha\rhoh)}
  \int_0^{1/y} u^{\alpha\rhoh-1} (1-u)^{\alpha\rho-1} \, \dd u \\
  &=& \frac{\Gamma(\alpha)}{\Gamma(\alpha\rho)\Gamma(\alpha\rhoh)}
    \int_y^\infty t^{-\alpha} (t-1)^{\alpha\rho-1} \, \dd t ,
  \yesnumber\label{e:law2}
\end{eqnarr*}
which gives us the density.

Now we substitute \eqref{e:law1} and \eqref{e:law2} into \eqref{e:inter1}:
\[ 
  \uo(x,y)
  = \frac{1}{\Gamma(\alpha\rho)\Gamma(\alpha\rhoh)}
  x^{\alpha\rhoh-1} y^{\alpha\rho}
  \int_{\frac{y}{x} \join 1}^{\frac{1}{x}}
  t^{-\alpha}
  (t-1)^{\alpha\rho-1}
  \biggl(t-\frac{y}{x}\biggr)^{\alpha\rhoh-1} \, \dd t . \]

The expression in the statement follows by a short manipulation
of this integral.
\end{proofof} 

\begin{proofof}{\autoref{Y (0,1)}}
According to \citet[Theorem 4.1]{Bau-reflected}, since $X$ is regular upwards,
we have the following formula for $r_1(0,y)$:
\[ r_1(0,y)
  = \lim_{z \downto 0}
  \frac{\uo(z,y)}{\stP_z(\tau_1^+ < \tau_0^-)} .
\]
We have found $\uo$ above, and as we already mentioned,
we have from \citet{Rog-hit} that
\begin{eqnarr*}
  \stP_x(\tau_1^+ < \tau_0^-)
  &=& \frac{\Gamma(\alpha)}{\Gamma(\alpha\rho)\Gamma(\alpha\rhoh)}
  \int_0^{x} t^{\alpha\rhoh-1}(1-t)^{\alpha\rho-1} \, \dd t .
\end{eqnarr*}
We may then make the following calculation, using l'H\^opital's rule on
the second line since the integrals converge,
\begin{eqnarr*}
  r_1(0,y)
  &=& \frac{1}{\Gamma(\alpha)} y^{\alpha-1}
  \lim_{z \downto 0}
  \frac{%
    \dint_0^{\frac{z(1-y)}{y-z}}
    s^{\alpha\rhoh-1} (s+1)^{\alpha\rho-1} \, \dd s
  }{%
    \dint_0^{z}
    t^{\alpha\rhoh-1}(1-t)^{\alpha\rho-1} \, \dd t
  } \\
  &=& \frac{1}{\Gamma(\alpha)} y^{\alpha-1}
  \lim_{z \downto 0}
  \frac{%
    z^{\alpha\rhoh-1} (1-y)^{\alpha\rhoh-1} (y-z)^{1-\alpha\rhoh}
    \frac{\partial}{\partial z}\bigl[ \frac{z(1-y)}{y-z} \bigr]
  }{%
    z^{\alpha\rhoh-1} \frac{\partial}{\partial z}\bigl[ z \bigr]
  } \\
  &=& \frac{1}{\Gamma(\alpha)} y^{\alpha\rho-1} (1-y)^{\alpha\rhoh}.
\end{eqnarr*}
Finally, the full potential density $r_1(x,y)$ follows simply by
substituting in the following formula, from the same theorem in
\cite{Bau-reflected}:
\[ r_1(x,y) = \uo(x,y) + \stP_x(\tau_0^- < \tau_1^+) r_1(0,y) .
\qedhere \]
\end{proofof} 

\bibliography{master}

\begin{thebibliography}{14}
\providecommand{\natexlab}[1]{#1}
\providecommand{\url}[1]{\texttt{#1}}
\expandafter\ifx\csname urlstyle\endcsname\relax
  \providecommand{\doi}[1]{doi: #1}\else
  \providecommand{\doi}{doi: \begingroup \urlstyle{rm}\Url}\fi

\bibitem[Baurdoux(2009)]{Bau-reflected}
E.~J. Baurdoux.
\newblock Some excursion calculations for reflected {L}\'evy processes.
\newblock \emph{ALEA Lat. Am. J. Probab. Math. Stat.}, 6:\penalty0 149--162,
  2009.
\newblock ISSN 1980-0436.

\bibitem[Bertoin(1996)]{Ber-Levy}
J.~Bertoin.
\newblock \emph{L\'evy processes}, volume 121 of \emph{Cambridge Tracts in
  Mathematics}.
\newblock Cambridge University Press, Cambridge, 1996.
\newblock ISBN 0-521-56243-0.

\bibitem[Blumenthal and Getoor(1968)]{BG-mppt}
R.~M. Blumenthal and R.~K. Getoor.
\newblock \emph{Markov processes and potential theory}.
\newblock Pure and Applied Mathematics, Vol. 29. Academic Press, New York,
  1968.

\bibitem[Blumenthal et~al.(1961)Blumenthal, Getoor, and Ray]{BGR-st-hit}
R.~M. Blumenthal, R.~K. Getoor, and D.~B. Ray.
\newblock On the distribution of first hits for the symmetric stable processes.
\newblock \emph{Trans. Amer. Math. Soc.}, 99:\penalty0 540--554, 1961.
\newblock ISSN 0002-9947.

\bibitem[Caballero and Chaumont(2006)]{CC-LS}
M.~E. Caballero and L.~Chaumont.
\newblock Conditioned stable {L}\'evy processes and the {L}amperti
  representation.
\newblock \emph{J. Appl. Probab.}, 43\penalty0 (4):\penalty0 967--983, 2006.
\newblock ISSN 0021-9002.
\newblock \doi{10.1239/jap/1165505201}.

\bibitem[Doney and Kyprianou(2006)]{DK-os}
R.~A. Doney and A.~E. Kyprianou.
\newblock Overshoots and undershoots of {L\'evy} processes.
\newblock \emph{The Annals of Applied Probability}, 16\penalty0 (1):\penalty0
  91--106, 2006.
\newblock \doi{10.1214/105051605000000647}.

\bibitem[Gradshteyn and Ryzhik(2007)]{GR}
I.~S. Gradshteyn and I.~M. Ryzhik.
\newblock \emph{Table of integrals, series, and products}.
\newblock Elsevier/Academic Press, Amsterdam, seventh edition, 2007.
\newblock ISBN 978-0-12-373637-6.
\newblock Translated from the Russian. Translation edited and with a preface by
  Alan Jeffrey and Daniel Zwillinger.

\bibitem[Kuznetsov and Pardo(2010)]{KP-HG}
A.~Kuznetsov and J.~C. Pardo.
\newblock Fluctuations of stable processes and exponential functionals of
  hypergeometric {L}\'evy processes.
\newblock \href{http://arxiv.org/abs/1012.0817v1}{arXiv:1012.0817v1} [math.PR],
  2010.

\bibitem[Kyprianou(2006{\natexlab{a}})]{Kyp}
A.~E. Kyprianou.
\newblock \emph{Introductory lectures on fluctuations of {L}\'evy processes
  with applications}.
\newblock Universitext. Springer-Verlag, Berlin, 2006{\natexlab{a}}.
\newblock ISBN 978-3-540-31342-7.

\bibitem[Kyprianou(2006{\natexlab{b}})]{Kyp-reflected}
A.~E. Kyprianou.
\newblock First passage of reflected strictly stable processes.
\newblock \emph{ALEA Lat. Am. J. Probab. Math. Stat.}, 2:\penalty0 119--123,
  2006{\natexlab{b}}.
\newblock ISSN 1980-0436.

\bibitem[Kyprianou et~al.(2010)Kyprianou, Pardo, and Rivero]{KPR-n-tuple}
A.~E. Kyprianou, J.~C. Pardo, and V.~Rivero.
\newblock Exact and asymptotic {$n$}-tuple laws at first and last passage.
\newblock \emph{Ann. Appl. Probab.}, 20\penalty0 (2):\penalty0 522--564, 2010.
\newblock ISSN 1050-5164.
\newblock \doi{10.1214/09-AAP626}.

\bibitem[Lamperti(1972)]{Lam-ssLT}
J.~Lamperti.
\newblock Semi-stable {M}arkov processes. {I}.
\newblock \emph{Z. Wahrscheinlichkeitstheorie und Verw. Gebiete}, 22:\penalty0
  205--225, 1972.

\bibitem[Rogozin(1972)]{Rog-hit}
B.~A. Rogozin.
\newblock The distribution of the first hit for stable and asymptotically
  stable walks on an interval.
\newblock \emph{Theory of Probability and its Applications}, 17\penalty0
  (2):\penalty0 332--338, 1972.
\newblock ISSN 0040585X.
\newblock \doi{10.1137/1117035}.

\bibitem[Vuolle-Apiala(1994)]{VA-Ito}
J.~Vuolle-Apiala.
\newblock It\^o excursion theory for self-similar {M}arkov processes.
\newblock \emph{Ann. Probab.}, 22\penalty0 (2):\penalty0 546--565, 1994.
\newblock ISSN 0091-1798.

\end{thebibliography}

\end{document}